\documentclass[a4paper,twoside,10pt]{amsart}
\usepackage{amsmath,amsthm,amssymb}
\usepackage{enumerate}
\usepackage{a4wide}
\usepackage{todonotes}

\usepackage{xcolor}
\usepackage[colorlinks=true, urlcolor=black, citecolor=black, linkcolor=black, hyperfootnotes=true]{hyperref}

\usepackage{url}

\usepackage[english]{babel}
\usepackage[utf8]{inputenc}
\usepackage[a4paper]{geometry}
\usepackage{latexsym}
\usepackage{mathtools}
\usepackage{amscd}
\usepackage{graphics}
\usepackage{color}
\usepackage{array}
\usepackage{mathrsfs}
\usepackage{graphicx}
\usepackage{stmaryrd}
\usepackage{comment}
\usepackage{varioref}
\usepackage{prettyref}

\newtheorem{thm}{Theorem}[section]
\newtheorem{cor}[thm]{Corollary}

\newtheorem{prop}[thm]{Proposition}

\numberwithin{claimcounter}{thm}

\newrefformat{sec}{Section \ref{#1}}
\newrefformat{thm}{Theorem \ref{#1}}
\newrefformat{cor}{Corollary \ref{#1}}
\newrefformat{lem}{Lemma \ref{#1}}
\newrefformat{prob}{Problem \ref{#1}}
\newrefformat{prop}{Proposition \ref{#1}}
\newrefformat{conj}{Conjecture \ref{#1}}
\newrefformat{fact}{Fact \ref{#1}}
\newrefformat{question}{Question \ref{#1}}

\theoremstyle{definition}
\newtheorem{defin}[thm]{Definition}

\numberwithin{equation}{section}

\newcommand{\N}{\mathbb{N}}
\newcommand{\R}{\mathbb{R}}

\newcommand{\Usph}{\mathbb{S}}
\newcommand{\comp}{\mathcal{K}}
\newcommand{\SCK}{\mathcal{SCK}}
\newcommand{\cP}{\mathcal{P}}

\DeclareMathOperator{\proj}{proj}

\DeclareMathOperator{\Emb}{Emb}

\DeclareMathOperator{\Ball}{B}

\DeclareMathOperator{\diam}{diam}
\DeclareMathOperator{\Pum}{Pum}
\DeclareMathOperator{\scc}{sc}
\renewcommand{\sc}{\scc}
\DeclareMathOperator{\PP}{PP}

\makeatletter
\newcommand{\alaligne}{~\vspace*{\topsep}\nobreak\@afterheading}
\makeatother

\usepackage{glossaries}

\newacronym{hi}{HI}{hereditarily indecomposable}

\usepackage[multiple]{footmisc}
\usepackage{bigfoot}

\DeclareNewFootnote{AAffil}[arabic]
\DeclareNewFootnote{ANote}[fnsymbol]

\begin{document}

%%%%% To ease editing, add:

%\baselineskip=17pt

%%%%%%%%%%%%%%%%

%% In the running head, give an abbreviation of the title. 
%\titlerunning{Ramsey theory without pigeonhole principle}

%\AtEndDocument{%
%  \par
%  \bigskip
%    \textsc{\footnotesize N. de Rancourt, Institut de Mathématiques de Jussieu - Paris Rive Gauche, Universit\'e Paris Diderot, Bo\^ite Courrier 7012, 75205 Paris Cedex 13, France}\\
%    \textit{E-mail address}: \texttt{derancou@dma.ens.fr}}

\title{Big Ramsey degrees in the metric setting}

\author[T. Bice, N. de Rancourt, J. Hubi\v{c}ka, M. Kone\v{c}n\'{y}]{T. Bice\textsuperscript{1}, N. de Rancourt\textsuperscript{1}, J. Hubi\v{c}ka\textsuperscript{1}, M. Kone\v{c}n\'{y}\textsuperscript{1}\\ \phantom{o}\\\textsuperscript{1}\emph{Institute of Mathematics of the Czech Academy of Sciences}\\ \phantom{o}\\\textsuperscript{2}\emph{University of Lille}\\ \phantom{o}\\\textsuperscript{3}\emph{Charles University, Prague}}

\date{}

\maketitle

%\address{}

%\subjclass{Primary: 46B20; Secondary: 46B03, 46B04, 46B28, 47A53.}

%\ackn{}

%\keywords{}

%%%%%%%%

\begin{abstract}
\emph{Oscillation stability} is an important concept in Banach space theory which happens to be closely connected to discrete Ramsey theory. For example, Gowers proved oscillation stability for the Banach space $c_0$ using his now famous Ramsey theorem for $\mathrm{FIN}_k$ as the key ingredient.
We develop the theory behind this connexion and introduce the notion of compact big Ramsey degrees, extending the theory of (discrete) big Ramsey degrees. We then prove existence of compact big Ramsey degrees for the Banach space $\ell_\infty$ and the Urysohn sphere, with an explicit characterization in the case of $\ell_\infty$.
\end{abstract}

\section{Introduction}
A \textit{discrete structure} is a structure in a relational language without unary predicates (e.g. orders, graphs, hypergraphs, etc.) Given two discrete structures $X$ and $Y$ in the same language we will denote by $X \choose Y$ the set of all embeddings $Y \to X$. Given a discrete structure $X$ and a finite substructure $A$ of $X$ we say that $A$ \textit{has a finite big Ramsey degree} in $X$ if there exists $t \geqslant 1$ such that every finite colouring of $X \choose A$ attains at most $t$ colors on $f[X] \choose A$ for some well-choosen $f \in {X \choose X}$. In this case, the \textit{big Ramsey degree} of $A$ in $X$ is the least such $t$. We say that $X$ \textit{has finite big Ramsey degrees} if every finite substructure of $X$ has a finite big Ramsey degree in $X$.

The infinite Ramsey theorem says that all big Ramsey degrees in $(\N, \leqslant)$ are equal to $1$. However, this is usually not the case. By Hjorth~\cite{HjorthOscillation}, no infinite homogeneous structure has all big Ramsey degrees equal to 1. However, they can still be finite. The task of identifying such structures was initiated by unpublished result by Laver who proved finiteness of big Ramsey degrees of $(\mathbb{Q}, <)$, quickly followed by their exact computation by Devlin~\cite{Devlin}. Since then there has been a lot of progress~\cite{LafSauVuk,NVT2008,DobrinenTriangleFree,ZuckerAgedEmbeddings,dobrinen2019ramsey,mavsulovic2020finite,HubickaPartialOrder,Hubicka2020uniform,BCHKNV,Balko2021exact,Balko2023}. For an overview, see Dobrinen's survey \cite{DobrinenSurvey}.

Recall that $c_0$ is the Banach space of all sequences $x \colon \N \to \R$ tending to zero at infinity, with the norm $||x||_\infty \coloneq \sup_{n \in \N} |x(n)|$, and that the \textit{unit sphere} of a Banach space $X$ is $S_X \coloneq \{x \in X \mid \|x\| = 1\}$. Also recall that the \textit{Urysohn sphere} $\Usph$ is the unique complete and separable metric space of diameter $1$ containing copies of all finite metric spaces of diameter at most $1$ to be \textit{ultrahomogeneous}, that is, such that every isometry between finite subsets of $\Usph$ extends to an onto isometry of $\Usph$. By a \textit{compactum}, we mean a compact metric space. %The two aforementioned results are the following.

A discrete structure is \emph{indivisible} if the big Ramsey degree of a vertex in $X$ is equal to $1$. Mostly motivated by the \textit{distortion problem} from Banach space theory, the two first indivisibility-like results for metric structures (often called \textit{oscillation stability results} by Banach space theorists), stated below, have been proved by Gowers~\cite{GowersLipschitz}, and Nguyen Van Th\'e and Sauer~\cite{NVTSauer}. Here, Lipschitz maps can be seen as continuous colourings, and compactness is the right metric analogue of finiteness.

\begin{thm}[Gowers \cite{GowersLipschitz}]\label{thm:Gowersc0}
Let $K$ be a compactum and $\chi \colon S_{c_0} \to K$ be a Lipschitz map. For every $\varepsilon > 0$, there exists a linear isometric copy $X$ of $c_0$ in $c_0$ such that $\diam(\chi(S_X)) \leqslant \varepsilon$.
\end{thm}

\begin{thm}[Nguyen Van Th\'e--Sauer \cite{NVTSauer}]\label{thm:NVTS}
Let $K$ be a compactum and $\chi \colon \Usph \to K$ be a Lipschitz map. For every $\varepsilon > 0$, there exists an isometric copy $X$ of $\Usph$ in $\Usph$ such that $\diam(\chi(X)) \leqslant \varepsilon$.
\end{thm}

The proof of Theorem \ref{thm:Gowersc0} is based on discrete approximations and was the reason why Gowers proved his now well-known $\mathrm{FIN}_k$ theorem, which is the main ingredient in the proof. The proof of Theorem \ref{thm:NVTS} is also using discrete approximations (following a combinatorial strategy which was proposed earlier by Nguyen Van Th\'e and Lopez-Abad~\cite{NVTLopezAbad}) combined with indivisibility results for metric spaces with finitely many distances.

The similarity between those results and indivisibility results makes it is natural to ask whether a suitable version of the notion of big Ramsey degrees for metric structures could be defined. This question was already addressed in the well-known Kechris--Pestov--Todorcevic paper \cite[\S{} 11(F)]{kpt}, where a definition was suggested. However, this definition is quite restrictive and fails to capture most interesting structures beyond the discrete ones. The goal of this work is to provide a more general notion, and to demonstrate its suitability on examples such as the Banach space $\ell_\infty$ and the Urysohn sphere. Our work on the Urysohn sphere builds on results on big Ramsey degrees of homogeneous structures with forbidden cycles announced at Eurocomb 2021 \cite{BCHKNV}.

Our motivations are twofold. First, Zucker~\cite{ZuckerCorrespondance} extended the KPT correspondence~\cite{kpt} to big Ramsey degrees, giving a correspondence between them and some dynamical invariants of automorphism groups. Our extension to the metric setting could allow us to study the same dynamical invariants for the automorphism groups of metric structures; no tool is currently available for studying those. Second, our methods could lead to a systematical study of the distortion phenomenon in Banach space theory, closely related to oscillation stability and not yet well understood. For instance, Odell and Schlumprecht's solution to the distortion problem~\cite{OdellSchlumprecht} show that the analogue of Theorem~\ref{thm:Gowersc0} fails for the separable Hilbert space. Metric big Ramsey degrees could help to express a quantitative and optimal version of their result.

%  Our long-term motivations are twofold.
% \begin{itemize}
%     \item The KPT-correspondence~\cite{kpt} between small Ramsey degrees of ultrahomogeneous discrete structures and dynamical properties of their automorphisms group has been extended by Zucker~\cite{ZuckerCorrespondance} to big Ramsey degrees. Our extension of big Ramsey degrees to the metric setting could allow us to study the dynamical objects introduced by Zucker in the case of automorphism groups of metric structures exhibiting some form of ultrahomogeneity; no tool is currently available for studying those.
%     \item The methods we develop could allow us to study more systematically the distortion phenomenon in Banach space theory, which is not yet well understood. For instance, Odell and Schlumprecht~\cite{OdellSchlumprecht} show that the analogue of Theorem~\ref{thm:Gowersc0} fails when replacing $c_0$ with the separable Hilbert space. The theory of metric big Ramsey degrees could allow one to express a quantitative and optimal version of their result.
% \end{itemize}

\section{Compact big Ramsey degrees}\label{sec:Def}

We first review, in a more general setting, some results on discrete big Ramsey degrees to motivate our definitions in the metric case. Our setting will be this of a monoid $M$ acting by injections on a set $X$. The action $M \curvearrowright X$ \textit{has a finite big Ramsey degree} if there exists $t \geqslant 1$ such that every colouring of $X$ with finitely colours takes at most $t$ values on a set of the form $p \cdot X$, $p \in M$. In this case, the \textit{big Ramsey degree} of the action is the least such $t$. Observe that if $Y$ is a discrete structure and $A \subseteq Y$ a finite substructure, then taking $M \coloneq {Y \choose Y}$ and $X \coloneq {Y \choose A}$ and considering the action by left-composition, we recover the classical notion of the big Ramsey degree of $A$ in $Y$. For $k \geqslant 1$, denote by $[k]$ the set $\{1, \ldots, k\}$.

\begin{defin}
Fix $M \curvearrowright X$ as above, and $k \geqslant 1$. Say that a colouring $\chi \colon X \to [k]$ is:
\begin{itemize}
    \item \textit{persistent} if for every $p \in M$, $\chi(p \cdot X) = [k]$;
    \item \textit{universal} if for every $l \geqslant 1$, every colouring $\psi \colon X \to [l]$ and every $p \in M$, there exists $q \in M$ and $f \colon [k] \to [l]$ such that $\psi\restriction_{pq\cdot X} = f \circ \chi\restriction_{pq\cdot X}$;
    \item a \textit{big Ramsey colouring} (or a \textit{canonical partition}, following \cite{LafSauVuk}) if it is both persistent and universal.
\end{itemize}
\end{defin}

The proof of the following fact is elementary.

\begin{prop}
Suppose that the action $M \curvearrowright X$ has a finite big Ramsey degree. Then it admits a big Ramsey colouring. Moreover, the number of colours of such a colouring is always equal to the big Ramsey degree of the action.
\end{prop}

Our metric setting will be this of a monoid $M$ acting by (non-necessarily onto) isometries on a complete metric space $X$. Inspired by Theorems \ref{thm:Gowersc0} and \ref{thm:NVTS}, we will define a \textit{colouring} of $X$ as a $1$-Lipschitz map $X \to K$, where $K$ is a compactum (the Lipschitz constant $1$ is here to ensure some rigidity). We will also allow some $\varepsilon$-approximation in our results. %Compared to the discrete setting, the situation is more complicated as the analogues of several notions equivalent to the existence of a finite big Ramsey degree in the discrete setting do not seem to be equivalent anymore, at least in our knowledge.
%
%\smallskip
%
The order on $\N$ will be ``replaced'' with the quasiordering between compacta defined as follows: $K \leqslant L$ if there exists a $1$-Lispchitz surjection $L \to K$. It is a classical fact that $K \leqslant L$ and $L \leqslant K$ if and only if $K$ and $L$ are isometric. If $\chi$ and $\psi$ are two maps defined on the same set and taking values in the same metric space, we will denote by $d_\infty(f, g)$ the supremum distance between $f$ and $g$.

\begin{defin}\label{def:BRD}
    Say that a compactum $K$ is:
    \begin{itemize}
        \item \textit{universal} with respect to the action $M \curvearrowright X$ if for every compactum $L$, every colouring $\psi \colon X \to L$, and every $\varepsilon > 0$, there exists $q \in M$, a colouring $\chi \colon X \to K$, and a $1$-Lipschitz map $f \colon K \to L$ such that $d_\infty(\psi\restriction_{q\cdot X}, f\circ \chi\restriction_{q\cdot X}) \leqslant \varepsilon$;
        \item the \textit{big Ramsey degree of the action $M \curvearrowright X$} if it is a $\leqslant$-least universal compactum.
    \end{itemize}
    Say that the action $M \curvearrowright X$ \textit{has a compact big Ramsey degree} if it admits a big Ramsey degree in the above sense.
\end{defin}

The big Ramsey degree of an action, if it exists, is obviously unique, up to isometry.

\begin{defin}\label{def:BRC}
    Say that a colouring $\chi \colon X \to K$ is:
    \begin{itemize}
        \item \textit{persistent} if for every $p \in M$, $\chi(p\cdot X)$ is dense in $K$;
        \item \textit{universal} if for every compactum $L$, every colouring $\psi \colon X \to L$, every $p \in M$ and every $\varepsilon > 0$, there exists $q \in M$ and a $1$-Lipschitz map $f \colon K \to L$ such that $d_\infty(\psi\restriction_{pq\cdot X}, f\circ \chi\restriction_{pq \cdot X}) \leqslant \varepsilon$;
        \item a \textit{big Ramsey colouring} if it is both persistent and universal.
    \end{itemize}
\end{defin}

\begin{prop}
Suppose that $\chi \colon X \to K$ is a big Ramsey colouring for the action $M \curvearrowright X$. Then $K$ is the big Ramsey degree of this action.
\end{prop}

\begin{prop}
Consider the following statements:
\begin{enumerate}[(1)]
    \item the action $M \curvearrowright X$ admits a universal compactum;
    \item the action $M \curvearrowright X$ has a compact big Ramsey degree;
    \item the action $M \curvearrowright X$ admits a universal clouring;
    \item the action $M \curvearrowright X$ admits a big Ramsey colouring.
\end{enumerate}
Then the following implications hold: (4) $\implies$ (3) $\implies$ (2) $\implies$ (1).
\end{prop}

While the analogues of the implications above are equivalent in the discrete setting, we do not know whether any of the reverse implications hold in the metric setting. The most relevant notion seems to be the existence of a big Ramsey colouring as, in the discrete setting, it is the closest to Zucker's condition for getting interesting dynamical consequences~\cite{ZuckerCorrespondance}. Also, in all metric examples for which we have been able to prove the existence of a universal compactum, we could also prove the existence of a big Ramsey colouring.

We end this section with mentioning that endowing discrete structures with the metric where any two distinct points at distance $1$, we can ``embed'' the classical discrete setting for big Ramsey degrees in our metric setting, making the discrete setting a particular case of the metric setting.

%The metric notions developed above are in fact generalization of their discrete counterparts: Let $M$ be a monoid acting by injections on a discrete set $X$. We can make $X$ into a metric space by deciding that any two distinct elements of $X$ are at distance $1$. The following result then links big Ramsey degrees of $X$ in the discrete and the metric sense.

\begin{comment}
\begin{prop}
    Under the above assumptions, the following conditions are equivalent.
    \begin{enumerate}[(1)]
        \item The action $M \curvearrowright X$ has a finite big Ramsey degree in the discrete sense.
        \item The action $M \curvearrowright X$ admits a universal compactum in the metric sense.
        \item The action $M \curvearrowright X$ admits a big Ramsey colouring in the metric sense.
    \end{enumerate}
    Moreover, if the above conditions hold, denote by $K$ the big Ramsey degree of $M \curvearrowright X$ in the metric sense. Then the metric space $K$ is metrically discrete (in the sense that for any two distinct $x, y \in K$, we have $d(x, y) = 1$) and the cardinality of $K$ is the big Ramsey degree of $M \curvearrowright X$ in the discrete sense.
\end{prop}*
\end{comment}

\section{Banach spaces}\label{sec:Banach}

In this section we study big Ramsey degrees of the spaces $\ell_p$ and $c_0$. Instead of colouring (linear isometric) embeddings of finite-dimensional subspaces into the whole space, we will equivalently colour finite tuples of elements of its unit sphere, which makes the presentation easier. Given a Banach space $X$ and $d \geqslant 1$, the set $(S_X)^d$ will be endowed with the supremum distance. We denote by $\Emb(X)$ the monoid of all linear isometric embeddings of $X$ into itself.

For $1 \leqslant p < \infty$ and a sequence $x \colon \N \to \R$, we let $\|x\|_p \coloneq \left(\sum_{n \in \N}|x(n)|^p\right)^{\frac{1}{p}}$, and $\|x\|_\infty \coloneq \sup_{n \in \N} |x(n)|$; and for $1\leqslant p \leqslant \infty$, we denote by $\ell_p$ the Banach space of all sequences $x \colon \N \to \R$ such that $\|x\|_p < \infty$, endowed with the norm $\|\cdot\|_p$. The space $c_0$ is a particular subspace of $\ell_\infty$.

Gowers' theorem~\ref{thm:Gowersc0} is equivalent to saying that the action $\Emb(c_0) \curvearrowright S_{c_0}$ admits a big Ramsey degree which is a singleton. However, the situation is different in higher arities. To see this, colour a pair $(x, y) \in S_{c_0}^2$ of disjointly supported vectors by the number of times their supports intertwine. Then every block-subspace of $c_0$ meets an infinite number of colours. While this ``colouring'' is neither Lipschitz nor compactum-valued, the idea can be developed to prove the following result.

\begin{thm}
    Let $d \geqslant 2$. Then the action $\Emb({c_0}) \curvearrowright (S_{c_0})^d$ does not admit a universal compactum.
\end{thm}

As mentioned in the introduction, Odell and Schlumprecht~\cite{OdellSchlumprecht} proved that the separable Hilbert space $\ell_2$ does not satisfy an analogue of Theorem~\ref{thm:Gowersc0}. In fact, their paper immediately implies a stronger conclusion.

\begin{thm}
    Let $d \geqslant 1$ and $1\leqslant p < \infty$. Then the action $\Emb({\ell_p}) \curvearrowright (S_{\ell_p})^d$ does not admit a universal compactum.
\end{thm}

Thus, our theory of big Ramsey degrees in its current form is not suitable for expressing a quantitative and optimal version of the Odell--Schlumprecht theorem, assuming it exists. A theory of noncompact big Ramsey degrees would be needed for this; we believe that such a theory could be developed, and keep it in mind for a future project.

We now turn to $\ell_\infty$. In the rest of this section, we fix $d\geqslant 1$ and consider the action $\Emb(\ell_\infty) \curvearrowright (S_{\ell_\infty})^d$. Classical arguments from Banach space theory show that, even when $d = 1$, this action does not admit a universal compactum. However, the proof involves a diagonal argument based on the Axiom of Choice. In such cases, imposing a definability restriction on colourings often allows one to get positive results (see e.g.~\cite{StevoBook}). The right topology is the \textit{weak-* topology} (here, we refer to the one we get when seeing $\ell_\infty$ as the dual of $\ell_1$). We can then define the notions of a \textit{definable big Ramsey degree} and a \textit{definable big Ramsey colouring} for the above action by considering, in Definitions~\ref{def:BRD} and~\ref{def:BRC}, only colourings that are Borel, or even Suslin--measurable, when $(S_{\ell_\infty})^d$ is endowed with the $d$-th power of the weak-* topology. All results proved in Section~\ref{sec:Def} remain valid for these definable notions, and it turns out that we can prove the existence of a definable big Ramsey colouring for the action $\Emb(\ell_\infty) \curvearrowright (S_{\ell_\infty})^d$. In order to state our result, preliminary definitions are needed. 

Put $\Ball_d \coloneq [-1, 1]^d$ and endow this set with the supremum metric. %Denote by $\proj_i \colon \Ball^d \to [-1, 1]$ the $i$-th coordinate functional.
The $d$ entries of a tuple $x \in (S_{\ell_\infty})^d$ will be denoted by $x_1, \ldots, x_d$. We use a functional notation for elements of $\ell_\infty$, so that for $n \in \N$ and $i \in [d]$, the $n$-th entry of the vector $x_i$ will be denoted by $x_i(n)$. We can finally let, for each $n \in \N$, $x(n)$ be the $d$-tuple $(x_1(n), \ldots, x_d(n))$; it is an element of $\Ball_d$. In this way, elements of $(S_{\ell_\infty})^d$ can be seen as maps $\N \to \Ball_d$.
If $X$ is a metric space, denote by $\comp(X)$ the set of all nonempty compact subsets of $X$, and endow it with the \textit{Hausdorff metric} $d_H$ defined by $d_H(K, L)\coloneq \max(\sup_{x \in K} d(x, L), \sup_{y \in L} d(y, K))$. Denote by $\SCK(\Ball_d)$ the set of all nonempty symmetric, convex and compact subsets of $\Ball_d$, and see it as a metric subspace of $\comp(\Ball_d)$. If $A \subseteq \Ball_d$, denote by $\sc(A)$ the symmetric convex hull of the set $A$.

\begin{defin}
A \textit{$d$-pumpkin} is a compact subset $\cP \subseteq \SCK(\Ball_d)$ such that $\{0\} \in \cP$, there exists $C \in \cP$ such that for all $i \in [d]$, $\proj_i(C) = [-1, 1]$, and the inclusion induces a dense linear order on $\cP$. We denote by $\Pum_d$ the set of all $d$-pumpkins, seen as a subset of $\comp(\SCK(\Ball_d))$.
\end{defin}

A $d$-pumpkin can be seen as a continuously growing symmetric compact convex subset of $\Ball_d$, starting at $\{0\}$ and such that the final step of the evolution touches all faces of the cube $\Ball_d$. It can be shown that the metric space $\Pum_d$ is compact.

\begin{defin}
For $x \in (S_{\ell_\infty})^d$, let:
$$\PP_d(x) := \Big\{\sc\{x(0), \ldots, x(n-1), tx(n)\} \,\Big|\, n \in \N, t \in [0, 1]\Big\} \cup \Big\{\overline{\sc\{x(n) \mid n \in \N\}}\Big\}.$$ This defines a definable colouring $\PP_d \coloneq (S_{\ell_\infty})^d \to \Pum_d$.
\end{defin}

This definition can be paraphrased as follows: the sets $\sc\{x(0), \ldots, x(n-1)\}$, $n \in \N$, must be steps of the evolution of the pumpkin $\PP_d(x)$, and the set $\overline{\sc\{x(n) \mid n \in \N\}}$ must be its final step. Between those steps, we ``fill in the holes'' in an affine way. %Note that this choice of affine extensions is totally arbitrary, and that other choices could have been done, as soon as they are coherent among all values of $x$. The big Ramsey colourings, when they exist, are far for being unique, after all.

\begin{thm}\label{thm:BRCellinfty}
    The colouring $\PP_d$ is a definable big Ramsey colouring of the action $\Emb(\ell_\infty) \curvearrowright (S_{\ell_\infty})^d$. In particular, $\Pum_d$ is the definable big Ramsey degree of this action.
\end{thm}

It is easy to see that $\Pum_1$ is a singleton. Thus, as a corollary of Theorem \ref{thm:BRCellinfty}, we get the following oscillation stability result for $\ell_\infty$, analogous to Theorem \ref{thm:Gowersc0}.

\begin{cor}\label{thm:OscStabellinfty}
Let $K$ be a compactum and $\chi \colon S_{\ell_\infty} \to K$ be a Lipschitz map that is also Borel (or Suslin-measurable) for the weak-* topology. Then for every $\varepsilon > 0$, there exists a linear isometric copy $X$ of $\ell_\infty$ in itself such that $\diam(\chi(S_X)) \leqslant \varepsilon$.
\end{cor}

\smallskip

The proof of Theorem \ref{thm:BRCellinfty} is based on the use of the Carlson--Simpson theorem. The natural presentation of $S_{\ell_\infty}$ as a set of infinite words over the alphabet $[-1, 1]$ makes its use particularly simple. Another ingredient in the proof is an analysis of the form of linear isometric copies of $\ell_\infty$ in itself, based on elementary Banach space theoretic tools.

\section{The Urysohn sphere}\label{sec:Urysohn}

Recall that $\Usph$ is the Urysohn sphere. As for Banach spaces, we will consider colourings of tuples from $\Usph$ rather than embeddings of finite substructures. For each $d \geqslant 1$, endow the $d$-th power $\Usph^d$ with the supremum metric. Denote by $\Emb(\Usph)$ the monoid of all (non-necessarily surjective) isometries of $\Usph$ into itself. The main result is the following.

\begin{thm}
    For every $d \geqslant 1$, the action $\Emb(\Usph) \curvearrowright \Usph^d$ admits a big Ramsey colouring.
\end{thm}

Our proof method is based on ideas developed in \cite{BCHKNV} for proving finiteness of the big Ramsey degrees of discrete versions of the Urysohn sphere. We don't work directly on $\Usph$ itself but on a metric space $\mathbb{T}$ that is bi-embeddable with it. This metric space is a well enough behaved space of sequences, allowing us the use of the Carlson--Simpson theorem. Our proof allows us to recover the fact that the big Ramsey degree of the action $\Emb(\Usph) \curvearrowright \Usph$ is a singleton, thus giving a new and short proof of Theorem \ref{thm:NVTS}, based on very different tools than the original proof. However, as soon as $d > 1$, part of our proof relies on a non-constructive argument, and we are currently not able to characterize the big Ramsey degrees completely. We are only able to give an upper bound of the big Ramsey degree in the sense of the quasiordering $\leqslant$, as a quotient of $\mathbb{T}^d$ by an action of the monoid of rigid surjections $\N \to \N$.

\smallskip

\bigskip\bigskip

\bibliographystyle{plain}
\bibliography{main}

\end{document}